\documentstyle[12pt]{article}
\newtheorem{thm}{Theorem}
\newtheorem{prop}{Proposition}

\newtheorem{lem}{Lemma}

\newtheorem{main}{Theorem}

\newenvironment{proof}{\medskip
\noindent {\bf Proof.}}{\hfill \rule{.5em}{1em}\mbox{}\bigskip}
\def\bea{\begin{eqnarray*}}
\def\eea{\end{eqnarray*}}
\def\bel{\begin{eqnarray}}
\def\eel{\end{eqnarray}}
\def\be{\begin{equation}}
\def\ee{\end{equation}}

\def\Bbb{\bf }
\title{Einstein Metrics on Complex Surfaces}
\author{Claude LeBrun
 \thanks{Supported
in part by  NSF grant DMS-92 .}\\SUNY Stony
 Brook}
\date{June, 1995}

\begin{document}
\maketitle

\section{Introduction}

Suppose  $M$ a compact manifold which admits an Einstein metric
 $g$ which is K\"ahler with respect to some
complex structure $J$. Is every other
Einstein metric $h$ on $M$ also K\"ahler-Einstein?
If the complex dimension of $(M,J)$ is
$\geq 3$, the answer is generally no;  for example,
${\Bbb CP}_3$  admits both the
 Fubini-Study metric, which is K\"ahler-Einstein,
and a non-K\"ahler  Einstein metric \cite{bes} obtained by appropriately
squashing the fibers of the twistor projection
${\Bbb CP}_3\to S^4$. Iterated Cartesian products with
${\Bbb CP}_1$ then provide counter-examples in all higher
dimensions.

However,
if $M$ is a 4-manifold,  so that
$(M,J)$ is a compact complex surface,
there is   reason to hope that the anwer
to the above question might    be
{\sl yes}.  Indeed,
 Hitchin \cite{hit2} was able to answer the question in the affirmative
for complex surfaces which admit {\em Ricci-flat}  K\"ahler  metrics;
his argument hinges on the fact that
 any
4-dimensional Einstein manifold satisfies
$$2\chi + 3\tau =\frac{1}{4\pi^2} \int \left(2|W_+|^2+\frac{s^2}{24}\right)
d\mu $$
where $s$ is the scalar curvature and $W_+$ is the
self-dual Weyl curvature, and on the observation that
   $|W_+|^2= s^2/24$ for any K\"ahler  surface.
 Much more recently, Seiberg-Witten theory \cite{KM,wit} has
provided new insights when
our K\"ahler-Einstein metric has $s <0$; in
this case the  K\"ahler-Einstein
metrics are   absolute minima of
the Riemannian functional $\int s^2 ~ d\mu$,
and a  close cousin of Hitchin's argument  therefore
 implies \cite{leb}
the desired result  for compact quotients $M={\Bbb C}{\cal H}_2/\Gamma$
of the unit ball in ${\Bbb C}^2$.

While   the answer to the above question regarding
Einstein 4-manifolds still  remains  elusive,
a related, narrower  problem is much more
tractable.  Namely, suppose that $(M,J)$ is a compact
complex surface with Hermitian metric $h$; that is, it is
supposed that  the Riemannian metric $h$ is $J$-invariant. If $h$ is an
Einstein
metric, is it necessarily K\"ahler with respect to $J$? In general, the answer
is
no; the {\em Page metric} \cite{page,bes} on
${\Bbb CP}_2\# \overline{\Bbb CP}_2$ is a counter-example. However,
as will be demonstrated in this note, this counter-example is
nearly unique:

\begin{main}
Let $(M^4, J)$ be a compact complex surface which admits an
Einstein metric $h$ which is Hermitian but not
  K\"ahler with
respect to $J$. Then $(M,J)$ is obtained from
${\Bbb CP}_2$ by blowing up one, two, or three points in general
position. Moreover, the isometry group of $h$ contains a
2-torus.
\end{main}
 In the one-point case, the proof will
also show that $(M,h)$ is precisely the Page metric, up to
isometry and rescaling.

The   proof of this result hinges upon the fact  that if
a Hermitian metric on a complex surface is
Einstein, it must be {\sl conformally}  K\"ahler; this
 follows from the combined results  of
Goldberg-Sachs \cite{gs} and Derdzinski \cite{der}.
Note that the anologous statement  is false for complex
manifolds of  higher dimension,  as is demonstrated by   the
    the ``squashed''
Einstein metric on ${\Bbb CP}_3$.

Specializing our initial question, we might now ask whether
  a compact complex surface $(M,J)$ can admit both a K\"ahler-Einstein metric
and an Einstein Hermitian metric which is {\em  not} K\"ahler.
The answer is no, unless perhaps if
$M = {\Bbb CP}_2\# 3 \overline{\Bbb CP}_2$.
 This follows
because one- and two-point point
blow-ups of ${\Bbb CP}_2$ have non-reductive
automorphism groups, and hence \cite{mat,bes} do not
carry  K\"ahler-Einstein metrics.
The  concluding section of this article will describe
a computional method of determining
whether every Einstein Hermitian metric
on ${\Bbb CP}_2\# 3\overline{\Bbb CP}_2$
is actually K\"ahler-Einstein.
 The same method may be applied to
the existence problem
for Einstein Hermitian metrics on
 ${\Bbb CP}_2\# 2\overline{\Bbb CP}_2$.

\noindent {\bf Acknowedgement.} The author
warmly thanks Jerzy Lewandowski for stimulating his
interest in the problem.

 \section{Einstein Hermitian Metrics}

In this section, we will study Einstein metrics which are
Hermitian with respect to some integrable complex structure
on a compact complex surface. These will, for the sake
of brevity, sometimes be referred to as Einstein
Hermitian metrics, so it is  worth warning the
 reader    that these
are  not {\em a priori}
Hermite-Einstein in the sense of the theory of
 holomorphic  vector bundles.

Let us begin with a local result concerning the
conformal curvature of Hermitian Einstein metrics:

\begin{lem}[Goldberg-Sachs] \label{gbs}
Let $(M,h)$ be an oriented
 Einstein 4-manifold. Assume that
there is an orientation-compatible
integrable complex structure $J$ on $M$ such that $h$ is
 $J$-invariant.
Then the self-dual Weyl curvature $W_+$ of $g$ is also
$J$-invariant. In particular, $W_+:\wedge^+\to \wedge^+$ has
at most 2 distinct eigenvalues at every point of $M$.
\end{lem}
A Lorentzian analogue of this result was first
discovered by Goldberg and Sachs \cite{gs}, but
two decades then elapsed  before it was realized \cite{pb,boy2}
that the same calculation concerning null involutive
sub-bundles of the complex
 tangent bundle
proves a  theorem
concerning  Riemannian signature metrics.
For a transparent spinorial proof, cf. \cite{pr} or \cite{nur}.

The self-dual Weyl curvature may be identified with
a symmetric trace-free  endomorphism of the
bundle $\wedge^+$ of self-dual 2-forms on our
Riemannian 4-manifold, and so has 3 eigenvalues at
each point.
 Under the action ${\Bbb Z}_2$-action
generated
an oriented orthogonal complex structure
$J$, however, the rank-3 bundle $\wedge^+$ decomposes into
irreducible sub-bundles of rank 1 and 2,
and two of the eigenvalues of a $J$-invariant
$W_+$ must therefore coincide.
The following
result of Derdzinski \cite[Theorem 2]{der}, however,
deals with Einstein manifolds with precisely this property.

\begin{lem}[Derdzinski] \label{derd}
Let $(M,h)$ be a connected oriented Einstein manifold such that
$W_+$ has at most 2 eigenvalues at each point. Then either
$W_+\equiv 0$, or else $W_+$ has exactly 2 distinct
eigenvalues at each point. In the latter case, moreover,
the conformally related metric $g=2\sqrt[3]{3}|W_+|^{2/3}h$ is locally
K\"ahler,
and is locally compatible with exactly one pair  $\pm J$ of  oriented
complex structures.
The scalar curvature $s$ of $g$ is then nowhere zero, and
$h=s^{-2}g$.
\end{lem}

The case in which
 our Hermitian metric $g$ satisfies $W_+\equiv 0$
may easily be handles by invoking
the work of Boyer:

\begin{lem} \label{bo}
Let $(M,J,h)$ be an Einstein Hermitian surface
with $W_+\equiv 0$.
Then $h$ is Ricci-flat and K\"ahler with respect to $J$.
\end{lem}
\begin{proof}
Since $(M,J,h)$ is Hermitian anti-self-dual, a result of Boyer \cite{boy1}
tells us that either $h$ is conformal to a scalar-flat K\"ahler metric,
or else that $b_1(M)=1$ and the conformal class $[g]$ has positive
Yamabe constant. The latter case, however, can be
excluded because our Einstein metric $h$ would have to
have positive scalar curvature, and hence positive
Ricci curvature; but this  would  imply that $b_1(M)=0$
by Bochner's theorem \cite{boc,bes}, and so lead to a contradiction. Hence $h$
is conformal to a scalar-flat K\"ahler metric, and, since it is the
Yamabe metric in its conformal class, it must therefore itself
be scalar-flat and K\"ahler.
\end{proof}

Combining these known facts now yields the following:

\begin{prop}
Let $(M,J,h)$ be a compact Einstein Hermitian manifold of
complex dimension 2. Then either $(M,J,h)$ is
a K\"ahler-Einstein manifold, or else there
 is an extremal  K\"ahler metric $g$ on
$(M,J)$ with  non-constant scalar curvature $s>0$ such that
$h=s^{-2}g$.
\end{prop}
\begin{proof} If  $W_+\equiv 0$,
Lemma \ref{bo} tells us that $h$ is K\"ahler, and we are done.
Otherwise,   Lemma \ref{derd} asserts that
the metric $g$ of is locally K\"ahler with respect to
exactly 2 complex structures, namely the two almost-complex structures
with respect to which $W_+$ is invariant; and by Lemma \ref{gbs},
the globally-defined
complex structure
$J$ is one of these.  Thus $(M,J,g)$ is a K\"ahler manifold.
But since $g$, being conformal to Einstein,
is a critical point of the conformally invariant functional
$\int |W_+|^2 d\mu= 6\pi^2 \tau  (M)+ \frac{1}{2} \int |W|^2 d\mu$,
and since $\int  |W_+|^2d\mu = \frac{1}{24}\int s^2d \mu$ for any
K\"ahler metric, it follows that $g$ is an extremal K\"ahler metric
in the sense of Calabi. Thus $\xi=J\mbox{grad}_gs$ is a Killing field of
$g$, and hence of $h=s^{-2}g$.  Now a result of Bochner \cite{boc,bes} says
that a compact manifold of non-positive Ricci curvature
can have a Killing field only if the field is parallel; but $\xi$  has
a zero at the minimum of $s^2$, and thus can be parallel only if
it is zero. Thus we either have $s=\mbox{const}$, in which case $h$
is K\"ahler-Einstein, or else the Einstein metric $h$ has
 positive scalar curvature. But if the latter happens,
the fact that $g$ is
in the same conformal class as $h$ implies that  its scalar curvature
$s\neq 0$  must also
be positive.
\end{proof}

While much of the above was already known to Derdzinski,
the next observation appears to be new:

  \begin{prop}
Let $(M^4, J)$ be a compact complex surface which admits an
Einstein metric $h$ which is Hermitian but not K\"ahler with respect to
$J$.   Then
the anti-canonical line bundle $K^{-1}$ of
$(M,J)$ is ample.
\end{prop}

\begin{proof}
Let $r$ denote the Ricci curvature  of the K\"ahler metric $g$,
and let  $\hat{r}=\frac{k}{4}h$
denote the Ricci curvature of the Einstein metric
$h=s^{-2}g$; here the constant $k=\hat{s}$ is the scalar curvature of
$h$.  The standard formula \cite{bes} for the effect of   a
conformal change $g\mapsto \psi^2g$ on the Ricci curvature  tells us that
$$\hat{r}_{ab}-r_{ab}=2\psi\nabla_a\nabla_b\psi^{-1}-(\psi\Delta \psi^{-1}
+ 3|d\log \psi|^2) g_{ab} ~ , $$
where the length of  1-forms is measured with respect to $g$; in
 our case, we therefore have
$$\frac{k}{4}{h}_{ab}=
r_{ab} + 2s^{-1}\nabla_a\nabla_b s - (s^{-1}\Delta s + 3 |d\log s|^2)g_{ab}$$
and hence
$$ s^{-1}\Delta s + 2|d\log s |^2 = \frac{s-ks^{-2}}{6} .$$
It follows that
  \be\label{ric}
r_{ab} + 2s^{-1}\nabla_a\nabla_b s = (\frac{2s+ks^{-2}}{12}+
|d\log s|^2)g_{ab}.\ee
Since both the K\"ahler metric $g$ and its Ricci curvature
$r$ are invariant under the action of $J$ on the tangent space, this implies,
in particular, that   the Hessian of $s$ is also $J$-invariant:
\be \label{inv}
\nabla_a\nabla_b s = {J_a}^c{J_b}^d\nabla_c\nabla_d s ~ .
\ee
(We remark in passing that (\ref{inv}) is exactly equivalent to
Calabi's extremal K\"ahler metric condition
$$   \nabla_{\bar{\mu}}
\nabla ^{\nu}s = 0 ~ . )$$

Now $2\pi c_1=2\pi c_1(K^{-1})$ is the de Rham class of
the Ricci form $\rho$ of our K\"ahler metric
$g$, and this 2-form is related to the Ricci curvature by
$$r_{ab}= \rho_{ac}{J_b}^c ~ . $$
Because the scalar curvature $s$ is a smooth
positive function on $M$, yet another
de Rham representative of $2\pi c_1 (K^{-1})$ is the
$(1,1)$-form
$\hat{\rho}$ defined by
\bea\hat{\rho}&=&
\rho +2 i \partial \bar{\partial} \log s \\&=&
\rho +  dJd\log s ~ .
\eea
But the $(1,1)$-form
  $\hat{\rho}$ is  `positive,'
in the sense that the symmetric tensor
field $q$ defined by
$$q_{ab}= \hat{\rho}_{ac}{J_b}^c ~ $$
is everywhere positive-definite.
Indeed,
\bea \hat{\rho}_{ab} &=& \rho_{ab} + 2\nabla_{[a}{J^c}_{b]}\nabla_c\log s\\
&=&\rho_{ab} - {J_b}^c \nabla_a\nabla_c \log s + {J_a}^c \nabla_b\nabla_c
\log s ~,\eea
so that
\bea
q_{ab} &=& \hat{\rho}_{ad}{J_b}^d\\&=&
 r_{ab} + \nabla_a\nabla_b \log s + {J_a}^c{J_b}^d \nabla_c\nabla_d \log s
\\&=&   r_{ab} + s^{-1}\nabla_a\nabla_b   s + s^{-1}{J_a}^c{J_b}^d
\nabla_c\nabla_d  s
\\&&  \hphantom{r_{ab}}
-(\nabla_a \log s) \nabla_b \log s-{J_a}^c{J_b}^d(\nabla_c \log s) \nabla_d
\log s
 ~  .
\eea
Substitution from (\ref{ric}) and (\ref{inv}) thus yields
\bea q_{ab} &=&   r_{ab} + 2s^{-1}\nabla_a\nabla_b   s
-(\nabla_a \log s) \nabla_b \log s-{J_a}^c{J_b}^d(\nabla_c \log s) \nabla_d
\log s\\
&=& \frac{2s+ks^{-2}}{12}g_{ab}+s^{-2}\left[  |d s|^2 g_{ab}
-( d s)_a  (d  s)_b- (J d s)_a  (Jd  s)_b \right] ~ ,
  \eea
 which is  manifestly  positive-definite
 because  $s$ and $k$ are both  positive.
Hence $c_1(K^{-1})$ is represented by the positive $(1,1)$-form
$\hat{\rho}/2\pi$,
and    the Kodaira embedding theorem \cite{gh} therefore
tells us that $K^{-1}$
is   ample. \end{proof}

This immediately implies the following:

  \begin{thm}
Let $(M^4, J)$ be a compact complex surface which admits an
Einstein metric $h$ which is Hermitian but not
  K\"ahler with
respect to $J$. Then $(M,J)$ is obtained from
${\Bbb CP}_2$ by blowing up one, two, or three points in general
position. Moreover, the isometry group of $h$ contains a
2-torus.
\end{thm}
 Here `general position' means that no two points coincide and
no three are collinear. After a projective-linear transformation
of ${\Bbb CP}_2$, we may therefore assume that our collection of
points is a subset of $\{ [1:0:0], [0:1:0], [0:0:1]\}.$
There are thus only 3 possible biholomorphism types for $(M,J)$.

\begin{proof}
Since the  anti-canonical line bundle  of $(M,J)$ is ample,
surface classification \cite{bpv,hit} tells us that  $(M,J)$ is
either ${\Bbb CP}_1\times {\Bbb CP}_1$ or else is
obtained from ${\Bbb CP}_2$ by blowing up $k$
distinct points in general position, $0\leq k \leq 8$.
However,  we also know that
$(M,J)$ carries an extremal K\"ahler
metric $g$ of non-constant scalar curvature,
so the Lie algebra of holomorphic
vector fields must be non-semi-simple
(and in particular non-trivial).
This eliminates ${\Bbb CP}_2$, ${\Bbb CP}_1\times
{\Bbb CP}_1$ and the $k$-point general-position blow-ups of ${\Bbb CP}_2$
for which  $4\leq k\leq 8$.

 Thus $(M,J)$ must be
obtained from ${\Bbb CP}_2$ by blowing up
1, 2, or 3 points in general position.
Choose homogeneous co\"ordinates on
${\Bbb CP}_2$ so that the points in question are elements of
$\{ [1:0:0], [0:1:0], [0:0:1]\}$, and observe that the
$U(1)\times U(1)$ action defined on ${\Bbb CP}_2$ by
$$\left[ \begin{array}{ccc}
e^{i\theta}&&\\&e^{i\phi}&\\&&1
\end{array}\right]$$
then lifts to the blow-up $M$; thus the  automorphism
group of $(M,J)$ contains the compact subgroup
$U(1)\times U(1)$.  But  \cite{cal} the identity component
of the isometry group of
an extremal K\"ahler metric
is a maximal compact subgroup of the identity component
of the complex automorphism group; and since the maximal
compact is unique up  to conjugation,  a suitable   change of
homogeneous co\"ordinates will   make the
extremal K\"ahler metric $g$   invariant under
above torus action.  Since any isometry of $g$
is also an isometry of $h=s^{-2}g$,  it follows that the
isometry group of  $h$ also contains $U(1)\times U(1)$.
\end{proof}

If $(M,g)$ is a compact oriented Einstein 4-manifold
with
 holonomy $SO(4)$ for which
$W_+$ has at most 2 eigenvalues
at each point, Lemma \ref{derd} thus implies that
 $M$ is diffeomorphic to
${\Bbb CP}_2\# \overline{\Bbb CP}_2$,
$({\Bbb CP}_2\# \overline{\Bbb CP}_2)/{\Bbb Z}_2$,
 ${\Bbb CP}_2\# 2 \overline{\Bbb CP}_2$,
or ${\Bbb CP}_2\# 3 \overline{\Bbb CP}_2$.

\section{Critical K\"ahler Classes}

In the last section, we saw that that a non-K\"ahler
Einstein Hermitian surface must be of the form
$(M,J,s^{-2}g)$, where $(M,J)$ is obtained from
${\Bbb CP}_2$ by blowing up 1, 2, or 3 points in
general position, and where $g$ is an extremal
K\"ahler metric of non-constant scalar
curvature $s > 0$. In the one-point case,
$h=s^{-2}g$ must be the Page metric, up to isometry and rescaling,
 because the
isometry group of the extremal K\"ahler metric
$g$ necessarily contains $U(2)$. In the other cases,
we can learn a bit more by asking which
K\"ahler class might contain such a metric $g$.

Let $[\omega ]$ denote the  K\"ahler   class of our
putative metric $g$. Since $g$ is extremal, there is an
open neighborhood of $[\omega ]\in H^{1,1}(M)=H^2(M)$
of classes
which are represented by extremal
K\"ahler metrics obtained as deformations of $g$.
On this open set, consider the
functional $\cal A$ which assigns to each cohomology class the
 integral $\int s^2 d\mu$ of the square of the
scalar curvature of the corresponding extremal
K\"ahler metric. Then $[\omega ]$ is a critical point
of this functional, since $\int s^2 d\mu = 24\int |W_+|^2d\mu$
for any K\"ahler metric, and the conformally
Einstein metric $g$ is a critical point of $\int |W_+|^2d\mu$,
considered as a functional on the space of all Riemannian
metrics.

This would be a useless observation were it not for
the fact that $\cal A$ has an invariant meaning.
Indeed, if $[\omega ]$ is the K\"ahler class of
an arbitrary extremal K\"ahler metric, we have
$${\cal A}= s_0^2\int  d\mu + \int (s-s_0)^2 d\mu=
32\pi^2\frac{(c_1\cdot [\omega ])^2}{[\omega ]^2}
-{\cal F}(\xi, [\omega ])
$$
where $s_0$ is the average value of the
scalar curvature,
$\cal F$ denotes the Futaki functional, and
$\xi = \mbox{grad}^{1,0} s$ is the extremal
vector field of the class $[\omega ]$. The latter,
moreover, may be determined up to conjugation
even \cite{fm} without knowing the extremal metric
explicitly. It is enough, in fact, to be able
to calculate the Futaki invariant explicitly, and this
has been done elsewhere \cite{ls1} for the blow-up
of ${\Bbb CP}_2$ at three points in general position.

Let $(M,J)$ be the blow-up of ${\Bbb CP}_2$ at the points
 $[1:0:0]$, $[0:1:0]$, and $[0:0:1]$.
The three blown-up points and the proper transforms of the lines
joining them  form a hexagon  of $(-1)$-curves in $M$:
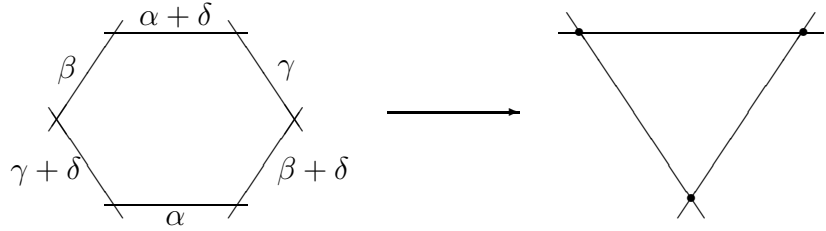
\begin{center}
\begin{picture}(240,80)(0,3)
\put(-7,70){\line(1,0){54}}
\put(40,75){\line(2,-3){28}}
\put(0,75){\line(-2,-3){28}}
\put(-7,5){\line(1,0){54}}
\put(40,0){\line(2,3){28}}
\put(0,0){\line(-2,3){28}}
\put(20,0){\makebox(0,0){$\alpha$}}
\put(20,77){\makebox(0,0){$\alpha+\delta$}}
\put(-21,56){\makebox(0,0){$\beta$}}
\put(62,56){\makebox(0,0){$\gamma$}}
\put(-29,18){\makebox(0,0){$\gamma+\delta$}}
\put(72,18){\makebox(0,0){$\beta+\delta$}}
\put(100,40){\vector(1,0){50}}
\put(210,0){\line(2,3){52}}
\put(220,0){\line(-2,3){52}}
\put(165,70){\line(1,0){100}}
\put(215,7.5){\circle*{3}}
\put(172.6,70.5){\circle*{3}}
\put(257.4,70.5){\circle*{3}}
\end{picture}
\end{center}
Since $b_2(M)=4$, there are two relations between these six curves---
namely, the three differences between opposite sides are homologous.
Thus, while the the areas $\alpha$, $\beta$ and $\gamma$ of the
three blow-up curves are independent, the only remaining
free parameter is the difference $\delta$ between the areas of
opposite sides of the hexagon. By performing a Cremona
transformation $$[z_0:z_1:z_2]\mapsto [1/z_0:1/z_1:1/z_2]$$
if necessary, we may arrange that $\delta \geq 0$, and
we will assume henceforth that this has been done.

Let us consider the hyperplane $P\subset
H^2({\Bbb CP}_2\# 3 \overline{\Bbb CP}_2)$
of K\"ahler classes
defined by the condition $\beta = \gamma$:
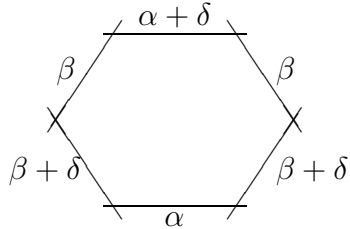
\begin{center}
\begin{picture}(100,80)(0,3)
\put(-7,70){\line(1,0){54}}
\put(40,75){\line(2,-3){28}}
\put(0,75){\line(-2,-3){28}}
\put(-7,5){\line(1,0){54}}
\put(40,0){\line(2,3){28}}
\put(0,0){\line(-2,3){28}}
\put(20,0){\makebox(0,0){$\alpha$}}
\put(20,77){\makebox(0,0){$\alpha+\delta$}}
\put(-21,56){\makebox(0,0){$\beta$}}
\put(62,56){\makebox(0,0){$\beta$}}
\put(-29,18){\makebox(0,0){$\beta+\delta$}}
\put(72,18){\makebox(0,0){$\beta+\delta$}}
\end{picture}
\end{center}
Since $\cal A$ is invariant under the
${\Bbb Z}_2$-action induced by
$$[z_0:z_1:z_2]\mapsto [z_0:z_2:z_1],$$
and because $P$ is exactly the fixed point
set of this action, any critical point
of ${\cal A}|_P$ is necessarily a
critical point
of ${\cal A}$, though the
converse of course need not be true.
Now it turns out that ${\cal A}|_P$
is rather easier to compute that ${\cal A}$,
and we shall therefore only consider this restricted
functional in the following discussion. It
should be emphasized, however, that this restriction
is completely {\em ad hoc}, and would have to be
eliminated in order to obtain a definitive treatment
of the problem.

For any K\"ahler class in $P$, the extremal K\"ahler  vector field $\xi$
must be invariant under $[z_0:z_1:z_2]\mapsto [z_0:z_2:z_1]$,
and so must be a multiple of the generator
 $\Xi$
of the ${\Bbb C}^\times$-action $[z_0:z_1:z_2]\mapsto [\zeta z_0: z_1: z_2]$.
Now $\Xi = \mbox{grad}^{1,0}t$ for a real-valued
function $t$ which, by the methods of \cite{ls1},
and preferably with the aid of a symbolic-manipulation program
such as {\sl Maple},
can be shown to satisfy
\def\a{\alpha}
\def\b{\beta}
\def\d{\delta}
\begin{eqnarray} [12\pi \omega ]^2 \int (t-t_0)^2 d\mu &= &
360 \b^3 \a \d^2 + 193\b^3\a^2\d + 276^4 \a\d  + 216 \b^2\d^3 \a
\nonumber
\\ &&
+60 \b\d^4\a + 48\b\d^3\a^2 +\d^6+12\b^6 + 96\b^4\a^2
{}~~~~~~~~
\label{rhs}
\\ && +72\b^5\a + 144\b^2\d^2\a^2 + 120\b^3\d^3
+ 138\b^4\d^2
\nonumber
\\ &&\nonumber +72\b^5\d
+ 54\b^2\d^4 + 12\b\d^5 + 6\d^5\a + 6\d^4\a^2
\end{eqnarray}
where $t_0$ is the average value of $t$.
On the other hand, it was shown in \cite{ls1} that
$$[\omega ]^2{\cal F}(\Xi, [\omega ])= 4(\beta-\alpha )\delta
\left [\frac{\delta^2}{3}+
\beta\delta+
\beta^2 \right] . $$
Since ${\cal F}(\Xi, [\omega ])=-\int (t-t_0) (s-s_0)d\mu$,
an explicit formula for ${\cal A}|_P$ can now be deduced by
setting $(s-s_0)= \lambda (t-t_0)$ and solving for
$\lambda$.
 The upshot is that $({\cal A}|_P)/{96\pi^2}$
is the quotient of two homogeneous sextics with integer coefficients, namely
\bea && 32\b^6 +160 \beta^5 \a+ 176 \b^5\d + 318 \b^4\d^2
+ 136\b^4\a^2+536\b^4\a\d   \\ && + 32\b^3\a^3 + 280\b^3\d^3
+
696\b^3\a\d^2 + 320\b^3\a^2\d  + 440\b^2\d^3\a \\ &&
+ 276\b^2\d^2\a^2 + 48 \b^2 \a^3\d   + 132 \b^2\d^4 + 32\b\d^5
+ 104 \b\d^3\a^2 \\ && + 24\b\a^3\d^2 + 136\b\d^4\a + 4\a^3\d^3 + 14\d^4\a^2
+ 16\d^5\a + 3\d^6
\eea
divided by the right-hand side of (\ref{rhs}).
Notice that
 ${\cal A}|_P$ in particular has homogeneity 0, as it must be
since $\int s^2d\mu$ is scale-invariant in  dimension 4.

To double-check our fomul{\ae}, let us first revisit the
case of ${\Bbb CP}_2\#\overline{\Bbb CP}_2$. This can be done
by simply setting $\beta =0$ in the above expressions, so that
$$\frac{1}{4\pi^2}\int_{{\Bbb CP}_2\#\overline{\Bbb CP}_2}
\frac{s^2}{24}d\mu =\frac{4 + 14x +16x^2 + 3x^3}{
x(6 + 6x + x^2)},
$$
where $x=\d/\a$. For $x > 0$, this has a unique critical point,
its absolute minimum, when $x= 2.183933404\ldots$
The Page metric must be conformal to an extremal
K\"ahler metric in this class,
for which  the area of a  projective line is
$3.183933404\ldots$
times that of the exceptional
divisor.  This agrees with the figure
obtained by more direct calculation;
cf.  \cite[p.338]{bes}.

Next we consider the 2-point blow-up
 ${\Bbb CP}_2\# 2 \overline{\Bbb CP}_2$ by instead
setting $\alpha = 0$ in our formula for
${\cal A}|_P$. This gives us
$$\frac{1}{4\pi^2}\int_{{\Bbb CP}_2\# 2 \overline{\Bbb CP}_2}
\frac{s^2}{24}d\mu =
\frac{32+ 176y+318y^2+280y^3+132y^4+32y^5+3y^6}{
12+72y+138y^2+120y^3+54y^4+12y^5+y^6}
$$
where $y=\d/\b$. For $y> 0$, this also has a unique critical point,
an absolute minimum, at $y=0.9577128052\dots$ and
for this critical K\"ahler class,
the area of a projective line is
$2.9577128052\dots$ times that of either exceptional divisor.
While it is unknown at present whether this K\"ahler
class is represented by an extremal metric, the trace-free part of the
Ricci tensor of such a metric would have to be rather small,
since one would then have
$$
\frac{1}{4\pi^2}\int_{{\Bbb CP}_2\#\overline{\Bbb CP}_2}
\left(2|W_+|^2+\frac{s^2}{24}\right)d\mu =
\frac{3}{4\pi^2}\int_{{\Bbb CP}_2\#\overline{\Bbb CP}_2}\frac{s^2}{24}d\mu =
7.136474469\ldots
$$
and the metric would thus satisfy
$$\frac{1}{8\pi^2}\int |\mbox{r}_0|^2d\mu=
0.136474469\ldots$$
by the Gauss-Bonnet formula for $2\chi + 3\tau =7$.
It would thus seem that there is at least a chance that
such a metric might in fact be conformal to Einstein.

Finally, we consider  ${\Bbb CP}_2\#3\overline{\Bbb CP}_2$.
In the region $\a,\b >0$, $\d \geq 0$, it appears
that ${\cal A}|_P$ has no critical points other than
an absolute minimum at $\a =\b$, $\d=0$,
which corresponds to multiples of the anti-canonical class.
Thus, at least if the symmetry condition $\b=\gamma$ is
imposed, it  seems that the only Einstein Hermitian
metrics on ${\Bbb CP}_2\#3\overline{\Bbb CP}_2$ are the
K\"ahler-Einstein metrics found by Siu \cite{s}; cf.
\cite{ty}. If this continues to hold even  when
$\b\neq \gamma$,  the uniqueness conjecture
of \S 1 will have survived an important test.

\pagebreak


\begin{thebibliography}{99}
 \bibitem{bpv}  W. Barth, C. Peters, and A. Van de Ven,
{\bf Compact Complex Surfaces}, Springer-Verlag, 1984.
\bibitem{bes} A. Besse, {\bf Einstein Manifolds}, Springer-Verlag, 1987.
\bibitem{boc} S. Bochner, {\em Vector Fields and
Ricci Curvature}, {\bf Bull.\ Am.\ Math. Soc.\ 52} (1946) 776--797.
\bibitem{boy1} C. Boyer, {\em
Conformal Duality and Compact Complex Surfaces},
{\bf Math.\ Ann.\ 274} (1986) 517-526.
\bibitem{boy2} C. Boyer, {\em Self-Dual and Anti-Self-Dual
Hermitian Einstein Metrics on Compact Complex Surfaces},
{\bf Contemp.\ Math.\ 71} (1988) 105--114.
\bibitem{cal} E. Calabi, {\em Extremal K\"ahler Metrics}, {\sl in}
{\bf Seminar
on Differential Geometry} (ed. S.-T. Yau), Princeton,  1982.
\bibitem{der} A. Derdzinski, {\em Self-Dual K\"ahler Manifolds and
Einstein Manifolds of Dimension Four}, {\bf Comp.\ Math.\ 49}
(1983) 405--433.
\bibitem{fm} A. Futaki and T. Mabuchi,
{\em Bilinear Forms and Extremal K\"ahler Vector Fields
Associated with K\"ahler Classes},
{\bf Math.\ Ann.\ 301} (1995) 199--210.
\bibitem{gh} P. Griffiths and J. Harris, {\bf Principles of
Algebraic Geometry},  Wiley-Interscience, 1978.
\bibitem{gs} J. Goldberg and R. Sachs, {\em
A Theorem on Petrov Types}, {\bf Acta Phys. Pol. {\sl Suppl.}
 22} (1962) 13--23.
\bibitem{hit}  N. Hitchin,  {\em On the Curvature of Rational
Surfaces}, {\bf Proc.\ Symposia  Pure Math.\  27}, pp 65-80,
1975.
\bibitem{hit2}  N. Hitchin,  {\em On Compact
Four-Dimensional Einstein Manifolds}, {\bf J. Diff.\ Geom.\ 9}
(1974) 435--442.
\bibitem{KM} P. Kronheimer, T. Mrowka,
{\em The Genus of Embedded Surfaces in the
Complex Projective Plane}, {\bf Math.\ Res.\ Letts.
1} (1994), 797--808.
\bibitem{mat}  Y. Matsushima, {\em Sur la Structure du
Groupe d'Hom\'eomorphismes
d'une Certaine Variet\'{e} K\"ahl\'erienne}, {\bf Nagoya Math. Journ.}
11 (1957) 145-150.
\bibitem{leb} C. LeBrun, {\em Einstein Metrics and Mostow Rigidity},
{\bf Math.\ Res.\ Letts. 2} (1995), 1--8.
\bibitem{ls1} C. LeBrun and
 S. Simanca, {\em Extremal K\"ahler Metrics and Complex Deformation Theory},
{Geom.\ Func.\ Analysis 4}(1994) 298--336.
\bibitem{ls} C. LeBrun and
 S. Simanca,
{\em On the K\"ahler Classes of Extremal  Metrics}, {\sl in}
{\bf Geometry and Global Analysis,}
{\sl Kotake, Nishikawa and Schoen, eds.}
pp. 255--271, Tohoku University, 1994.
\bibitem{nur} P. Nurrowski, {\bf Einstein Equations
and Cauchy-Riemann Geometry}, SISSA Ph.D. Thesis.
\bibitem{page} D. Page, {\em A Compact Rotating Gravitational
Instanton}, {\bf Phys.\ Lett.\ 79B} (1979) 235--238.
\bibitem{pb} M. Przanowski and B. Broda, {\em
Locally K\"ahler Gravitational Instantons},
{\bf Acta Phys. Pol. B 14} (1983) 637--661. 
\bibitem{pr} R. Penrose and W. Rindler,
{\bf Spinors and Space-Time}, vol. 2, Cambridge University Press, 1986.
\bibitem{s} Y.T. Siu, {\it The Existence of K\"ahler-Einstein Metrics on
Manifolds with Positive Anti-Canonical Line Bundle and Suitable Finite
Symmetry Group},
 {\bf Ann.\ Math.\ 127} (1988) 585-627.
\bibitem{ty} G. Tian and S.T. Yau,  {\it K\"ahler-Einstein Metrics on
Complex Surfaces with  ${\bf c}_1 > 0$},
 {\bf Comm.\ Math.\ Phys.\ 112} (1987) 175-203.
 \bibitem{wit} E. Witten, {\em Monopoles and Four-Manifolds},
{\bf Math.\ Res.\ Letts.
1} (1994), 769--796.


\end{thebibliography}
\end{document}